\documentclass[preprint,10pt]{article}
\usepackage[centertags]{amsmath}
\usepackage[colorlinks,linkcolor=red,anchorcolor=blue,citecolor=green]{hyperref}
\usepackage{amsfonts}
\usepackage{amsfonts}
\usepackage{newlfont}
\usepackage{graphicx}
\usepackage{subfigure}
\usepackage{stfloats}
\usepackage{authblk}

\numberwithin{equation}{section}
\newtheorem{theorem}{Theorem}[section]
\newtheorem{lemma}{Lemma}[section]

\newtheorem{example}[theorem]{Example}

\newcommand{\abs}[1]{\left\vert#1\right\vert}

\usepackage{epstopdf}
\usepackage{color}
\definecolor{huxl}{rgb}{0,0,0}   
\usepackage{lineno}

\begin{document}

\title{On a novel numerical quadrature based on cycle index of symmetric group for the Hadamard finite-part integrals}

\author{Jiajie Yao\thanks{Email:jjyao@zjut.edu.cn} }
\author{Congcong Xie\thanks{Corresponding author:ccxie@zjut.edu.cn} }

\affil{Department of Mathematics, Zhejiang University of Technology, $\textstyle{310023}$ Hangzhou, China}


\date{\today}

\maketitle

\begin{abstract}
{\color{huxl} 
	To evaluate the Hadamard finite-part integrals accurately, a novel interpolatory-type quadrature  is proposed in this article. 
	In our approach,  numerical divided difference is utilized to represent the high order derivatives of the integrated function, which make it possible to reduced the numerical quadrature into a  concise formula based on the cycle index for  symmetric group.}  In addition, convergence analysis is presented and the error estimation is given.  {\color{huxl} Numerical results are presented on cases with different weight functions, which substantiate the performance of the proposed method.
}
\\[3mm]
\textbf{keyword}:
Cauchy principal value integral; Hadamard finite-part integral; numerical divided difference; the cycle index of symmetric group.
\end{abstract}



\section{Introduction}
The numerical evaluation of Cauchy principal value integrals and Hadamard
finite-part integrals has received considerable attention,  especially in the
boundary element methods\cite{NM00,OW08,yu02}, where the efficiency of numerical evaluation of such Hadamard finite-part integrals are essential for the boundary element methods.
{\color{huxl} Considering} the numerical evaluation of the Hadamard finte-part integrals given by
\begin{eqnarray}\label{s1.1}
H_p(\omega;f;\xi)=\int_a^b\!\!\!\!\!\!\!\!\!=\omega(x)\frac{f(x)}{(x-\xi)^{p+1}}\mathrm{d}x,\
\ \xi\in (a,b),\ p\in \mathbb{N}_0:=\{0,1,\cdots\},
\end{eqnarray}
where $f(x)$ are assumed to be smooth functions and $\omega(x)$ are nonnegative weight functions on $[a,b]$.
In the case of $p=0$, the integral (\ref{s1.1}) {\color{huxl} is reduced to} the well-known Cauchy
principal value integral. {\color{huxl} In this sense}, Hadamard finite-part integrals {\color{huxl} could be also} considered as the {\color{huxl} generalization} of Cauchy principal value integrals.

A fundamental idea is to approximate the function $f(x)$ in (\ref{s1.1}) directly, such as the Lagrange or Lagrange-Hermite interpolation  based on a set of distinct nodes. There have been numerous studies on numerical evaluation of integrals with such singularities, such as Gaussian quadrature rule\cite{JCAM95,MC85,Milo2006,Mone97,Wang2017},
(composite) Newton-Cotes method\cite{JCAM10,JCAM101,Liu2017,NM08,JCAM091} and some other improved methods\cite{JCAM97,MC91,JCAM04,JCAM001,Keller2016,Setia2014,Xiang2016}. From a theoretical point of view,  the above mentioned methods turn out to be convergent only if the function $f$ is sufficiently smooth.  

The main difficulty arising in practical numerical evaluations is the so called "numerical cancellation", which  constantly happen  when the quadrature nodes approaching the singularity $\xi$.  Considering  numerical quadratures, as well as  mentioned in \cite{DR84,JCAM09,MC88}, subtracting out the singularity from $f(x)$,   (\ref{s1.1}) leads to
\begin{eqnarray}\label{s1.2}
\displaystyle H_p(\omega;f;\xi)=
\int^b_{a}\frac{\omega(x)}{(x-\xi)^{p+1}}\Big[f(x)-\sum_{j=0}^p\frac{f^{(j)}(\xi)}{j!}(x-\xi)^j\Big]\mathrm{d}x
+\displaystyle\sum_{j=0}^p\frac{f^{(j)}(\xi)}{j!}\int_a^b\!\!\!\!\!\!\!\!\!=\frac{\omega(x)}{(x-\xi)^{p+1-j}}\mathrm{d}x.
\end{eqnarray}
When considering the calculation of above integral, the second part of the right hand side are the Hadamard finite-part integrals, which  could be computed analytically with the weight function $\omega(x)$ known. The first part is a Riemann integral, which could be approximate it by general numerical quadrature rules, such as the Gaussian rule.  
However, the theoretical convergence rate by applying general numerical quadrature always breaks down since $\xi$ might be very close to certain quadrature nodes, which may lead to large rounding errors due to numerical cancellation (see \cite{JCAM95}). In this case, strong numerical cancellation could be presented and the quadrature rules might not be uniformly convergent for all $\xi\in (a,b)$. 

Some improved numerical strategies emerged recently to reduce the numerical cancellation of the approximation to (\ref{s1.2}) .  P. Kim and B.I.Yun \cite{JCAM001,JCAM002} construct a quadrature rule of
interpolatory-type based on the trigonometric interpolation for Cauchy principal value integrals.
The proposed rule is numerically stable and also estimates the error bounds.
Other authors for example, T. Hasegawa and T. Torii\cite{MC91} give an automatic
quadrature rule for computing Cauchy integrals.
They approximate the function $f$ by a sum of Chebyshev polynomials
whose coefficients are computed using the FFT. Similarly, later in \cite{JCAM04} and \cite{Hasa2017},
T. Hasegawa presents numerically stable interpolatory integration rules
to approximate Hadamard finite-part integrals and Cauchy principal value integrals with logarithmic singularity, respectively. All these results are extensions of the Clewshaw-Curtis quadrature rule. From another point of view, we can also interpolate function $f(x)$ on the nodes which are far from the singularity $\xi$ (see \cite{JCAM97}).
The method avoids the numerical cancellation by choosing better interpolation function but reduces the high precision of Gaussian quadrature rule.


The novelty of the current work are in two fold. Firstly, lie in the application of numerical divided difference,say $\dfrac{f^{(j)}(\xi)}{j!}=f[\underbrace{\xi,\cdots,\xi}_{j+1}]:=f[\xi^{j+1}] \quad \forall j=1,2,
\cdots$,  it  leads to a more concise formulation on calculation the first part of \eqref{s1.2}. It leads to
\begin{equation*}
 f(x)-\sum_{j=0}^p\frac{f^{(j)}(\xi)}{j!}(x-\xi)^j
	 =   f(x)-\sum_{j=0}^pf[\underbrace{\xi,\cdots,\xi}_{j+1}](x-\xi)^j 
	 =f[x,\xi^{p+1}](x-\xi)^{p+1},
\end{equation*}
where $f[x,\xi^{p+1}]$ denotes the divided difference of the function $f$ at the points $x,\xi$, 
and $\xi$ is repeated $p+1$ times. Hence, the equation(\ref{s1.2}) can be written as the following form
\begin{equation}\label{s1.3}
H_p(\omega;f;\xi)=
\int^b_{a}\omega(x)f[x,\xi^{p+1}]\mathrm{d}x+\sum_{j=0}^p\frac{f^{(j)}(\xi)}{j!}\int_a^b\!\!\!\!\!\!\!\!\!=\frac{\omega(x)}{(x-\xi)^{p+1-j}}\mathrm{d}x,
\end{equation}
where general Gaussian quadrature rule could be applied.  Secondly, the proposed quadrature rule can significantly alleviate the difficulty caused by numerical cancellation. Since Lagrange interpolation are calculated based on a different group of nodes, the distances between the nodes $a_0,a_1,\cdots,a_n$ are large enough so that the computation of the divided differences does not present numerical cancellation. 

The rest of the paper is organized as follows.  In order to get a explicit result, we need a conception of the cycle index of symmetric group, which is introduced in Section 2. A new quadrature rule for Hadamard finite-part integrals using the numerical divided difference rule is presented in section 3 while a convergent error estimation formula is given in section 4. Finally, we show some numerical examples to illustrate the effectiveness and accuracy in Section 5.

\section{ Cycle index for symmetric group}\label{sec:2}

Let $\mathfrak{G}_n$ be a symmetric group of degree $n$ ($n\in \mathbb{N}=\{1,2,\cdots\}$),
and for any permutation $\sigma\in \mathfrak{G}_n$, assume that $c_i(\sigma)$ is the number of the cycles of length $i$ in $\sigma$.
Now, we define the cycle index for symmetric group $\mathfrak{G}_n$ as follows:
\begin{eqnarray*}
Z_n(x_1,x_2,\cdots,x_n):=\frac{1}{|\mathfrak{G}_n|}\sum_{\sigma\in\mathfrak{G}_n}x_1^{c_1(\sigma)}x_2^{c_2(\sigma)}\cdots x_n^{c_n(\sigma)},
\end{eqnarray*}
where $|\mathfrak{G}_n|$ represents the order of $\mathfrak{G}_n$. For convenience, we denote it as
\begin{eqnarray*}
Z_n(x_i):=Z_n(x_i|1\leq i\leq n):=Z_n(x_1,x_2,\cdots,x_n).
\end{eqnarray*}

As we know $|\mathfrak{G}_n|=n!$, and the cycle index of the symmetric group can be expressed as the following explicit formulation\cite{NY72}:
\begin{eqnarray}\label{s2.1}
Z_n(x_i)=\sum_{a\in\pi_n}\frac{1}{1^{a_1}a_1!2^{a_2}a_2!\cdots n^{a_n}a_n!}x_1^{a_1}x_2^{a_2}\cdots x_n^{a_n},
\end{eqnarray}
where
\begin{eqnarray*}
\pi_n=\{(a_1,a_2,\cdots,a_n)\in \mathbb{N}_0\ |\
1a_1+2a_2+\cdots+na_n=n\}.
\end{eqnarray*}

In fact, the cycle index of the symmetric group $Z_n(x_i)$ is a polynomial of degree $n$ with the variable $x_1,x_2,\cdots,x_n$ .
However, it is difficult to obtain the expression from \eqref{s2.1} for large $n$, directly.
The following recursion relation for $Z_n(x_i)$ can be easily verified (cf. \cite{Riordan,SC07,Wang2006}):
\begin{eqnarray*}\label{s2.2}
nZ_n(x_i)=\sum_{j=1}^nx_jZ_{n-j}(x_i),\ \ \ n\geq 1;\ \ \ \ Z_0=1.
\end{eqnarray*}
Therefore, we have
\begin{equation*}
\begin{aligned}
&Z_1(x_1)=x_1,\\
&Z_2(x_i)=\frac{1}{2}(x_1^2+x_2),\\
&Z_3(x_i)=\frac{1}{6}(x_1^3+3x_1x_2+2x_3),\\
&Z_4(x_i)=\frac{1}{24}(x_1^4+6x_1^2x_2+3x_2^2+8x_1x_3+6x_4).
\end{aligned}
\end{equation*}

Now, let $\displaystyle\Omega(x)=\prod_{i=0}^n(x-a_i)$, $g(x)=\log|\Omega(x)|$. Then we have
\begin{equation}\label {s2.4}
g^{(k)}(x)=(-1)^{k-1}(k-1)!\sum_{i=0}^n(x-a_i)^{-k},\ \ \ x\neq a_i (i=0,1,\cdots,n).
\end{equation}
This is together with Fa$\grave{\mathrm{a}}$ di Bruno's formula
\begin{equation*}
\frac{\mathrm{d}^k}{\mathrm{d}x^k}\exp(g(x))=k!\exp(g(x))Z_k\Big(\frac{g^{(r)}(x)}{(r-1)!}\Big).
\end{equation*}

\begin{lemma}\label{lem2.1}
Assume that $a_0,a_1,\cdots,a_n$ are distinct points, then we have
\begin{eqnarray}\label{s2.3}
\frac{1}{k!}\Omega^{(k)}(x)=\Omega(x)Z_k(-S_r(x)),\ \ S_r(x)=\sum_{i=0}^n\frac{1}{(a_i-x)^r}\ \ x\neq a_i.
\end{eqnarray}
Especially,
\begin{eqnarray*}
\frac{1}{k!}\Omega^{(k)}(x)\Big |_{x=a_i}=\Omega'(a_i)Z_{k-1}(-S_{ri}(a_i)),\ \ S_{ri}(x)=\sum_{j=0\atop j\neq i}^n\frac{1}{(a_j-x)^r}.
\end{eqnarray*}
\end{lemma}
\textbf{Proof}. For $x\neq a_i$,
\begin{equation*}
\Omega^{(k)}(x)=\frac{\mathrm{d}^k}{\mathrm{d}x^k}(\mathrm{sgn}\Omega(x)\exp(g(x)))=k!\Omega(x)Z_k\Big(\frac{g^{(r)}(x)}{(r-1)!}\Big),
\end{equation*}
together with (\ref{s2.4}), we obtain (\ref{s2.3}).

For $x=a_i$, according to the Leibnitz formula for high order derivatives,
\begin{equation*}
\begin{aligned}
\Omega^{(k)}(x)&=(x-a_i)\Big(\frac{\Omega(x)}{x-a_i}\Big)^{(k)}+k\Big(\frac{\Omega(x)}{x-a_i}\Big)^{(k-1)},\\
&=k!\Omega(x)Z_k(-S_{ri}(x))+k!\frac{\Omega(x)}{x-a_i}Z_{k-1}(-S_{ri}(x)).
\end{aligned}
\end{equation*}
Applying $\Omega(a_i)=0$ and $\displaystyle\frac{\Omega(x)}{x-a_i}\Big|_{x=a_i}=\Omega'(a_i)$, we complete the proof. \hfill $\square$
\\

The next theorem is verified from Lemma \ref{lem2.1}.

\begin{theorem}\label{lem2.2}
We denote by $l_i(x)$ the Lagrange fundamental polynomial interpolating at the points $a_0,a_1,\cdots,a_n$,
\begin{eqnarray*}
l_i(x)=\frac{\Omega(x)}{(x-a_i)\Omega'(a_i)}.
\end{eqnarray*}
For $x\neq a_j\ (j\neq i)$, we have
\begin{eqnarray}\label{s2.5}
\frac{1}{k!}l^{(k)}_i(x)=l_i(x)Z_k(-S_{ri}(x)).
\end{eqnarray}
In particular for $x=a_j\ (j\neq i)$, the following equation holds,
\begin{eqnarray}\label{s2.7}
\frac{1}{k!}l^{(k)}_i(x)=l_i'(a_j)Z_{k-1}(-S_{rij}(x)),\ \ S_{rij}(x)=\sum_{l=0\atop l\neq i,j}^n\frac{1}{(a_l-x)^r}.
\end{eqnarray}
\end{theorem}
\textbf{Proof}.
Due to $\displaystyle l_i(x)=\frac{1}{\Omega'(a_i)}\prod_{j=0\atop j\neq i}^n(x-a_j)$, we obtain (\ref{s2.5}) combined with (\ref{s2.3}).

For $x=a_j\ (j\neq i)$, following from the Leibnitz formula,
\begin{eqnarray*}
\begin{aligned}
l_i^{(k)}(x)&=(x-a_j)\Big(\frac{l_i(x)}{(x-a_j)}\Big)^{(k)}+k\Big(\frac{l_i(x)}{(x-a_j)}\Big)^{(k-1)}\\
&=k!l_i(x)Z_k\Big(-\sum_{l=0\atop l\neq i,j}^n\frac{1}{(a_l-x)^r}\Big)+k!\frac{l_i(x)}{(x-a_j)}Z_{k-1}\Big(-\sum_{l=0\atop l\neq i,j}^n\frac{1}{(a_l-x)^r}\Big).
\end{aligned}
\end{eqnarray*}
Applying $l_i(a_j)=0\ (j\neq i)$, \eqref{s2.7} holds and the proof is completed. \hfill $\square$


\section{Quadrature formula}


Let $x_1,x_2,\cdots,x_m$ be quadrature nodes on $[a,b]$, which is  based on the zeros of the orthonormal polynomial with respect to the weight $\omega(x)$.
 Gaussian quadrature rule is applied to the fist integral on the right hand of \eqref{s1.3}, which leads to
\begin{eqnarray}\label{s3.1}
\displaystyle H_m(\omega;f;\xi)=
\sum_{k=1}^m\lambda_kf[x_k,\xi^{p+1}]+\sum_{j=0}^p\frac{f^{(j)}(\xi)}{j!}\int_a^b\!\!\!\!\!\!\!\!\!=\frac{\omega(x)}{(x-\xi)^{p+1-j}}\mathrm{d}x.
\end{eqnarray}
The divided differences $f[x_k,\xi^{p+1}]$ are obtained making use of
the value of function $f$ at $x_k$ and $\xi$. However, when $x_k$ is
very close to the singularity $\xi$, the practical evaluation
will produce large rounding errors due to numerical cancellation.
Therefore, it is expected that the distance between $x_k$ and $\xi$ is
large enough. But how to define "very close" and "large enough"?  For example, let $x_c$ be the quadrature node closest to $\xi$, then the distances of the
rest quadrature nodes to $\xi$ are considered to be large enough. That is, for any fixed $\xi$, $x_c$ is defined by
\begin{eqnarray*}
|\xi-x_c|=\min\{|\xi-x_k|, k=1,2,\ldots,m\},
\end{eqnarray*}
when $\xi$ is equidistant between two zeros, i.e., $\xi=\frac{1}{2}(x_k+x_{k+1})$ for some $k\in\{1,2,\ldots,m-1\}$,
then we choose two nodes $x_c=x_k$ and $x_c=x_{k+1}$.

In this case, there is an effective way to calculate $f[x_k,\xi^{p+1}]$ proposed by X. H. Wang, H. Y. Wang and M. J. Lai\cite{SC05}, where  the divided differences of the interpolation of $f$ is used to approximate the divided differences of $f$ .  The divided differences $f[x_k,\xi^{p+1}]$ can be calculated
easily without producing numerical cancellation expect $f[x_c,\xi^{p+1}]$.
For $f[x_c,\xi^{p+1}]$ we can use the divided differences of the
 interpolation of $f$ on another group of nodes $a_0,a_1,\cdots,a_n$ and the
 distances between the nodes are large enough. We assume that
 $a_0,a_1,\cdots,a_n$ are quasi-uniform nodes, that is, there exists constant $c>0$ such that
\begin{eqnarray*}
\max\{\overline{a}_{i+1}-\overline{a}_{i}|\ 0\leq i<n\}\leq cd.
\end{eqnarray*}
Here $\overline{a}_0<\overline{a}_{1}<\cdots<\overline{a}_{n}$
is the rearrangement of  $a_0,a_1,\cdots,a_n$ and $d$ is defined by
\begin{eqnarray*}
d=\min\{\overline{a}_{i+1}-\overline{a}_{i}|\ 0\leq i<n\},
\end{eqnarray*}
obviously, we need $x_c,\xi\in
[\overline{a}_0,\overline{a}_n]$ and $d>|x_c-\xi|$.

Furthermore, how to choose the nodes $a_0,a_1,\cdots,a_n$ is also an important problem.
Actually, in numerical computing practice,
we choose the equidistant nodes in the interval $[a,b]$ and let $\xi$ be one of the nodes.
For instance $b-\xi< \xi-a$, then we suppose that the subinterval $[\xi,b]$ is divided into $\nu+1$ parts
and let $h=\frac{b-\xi}{\nu+1}$ be step.
In the subinterval $[a,\xi]$, it is divided into $n-\nu$ parts under the step $h$. Let
\begin{equation}\label{s3.2}
\begin{aligned}
a_0=\xi,&\ a_1=\xi+h,\ a_3=\xi+2h,\ \cdots,\
a_{2\nu-1}=\xi+\nu h,\\
&\ a_2=\xi-h,\ a_4=\xi-2h,\ \cdots,\
a_{2\nu}=\xi-\nu h,\\
&\ a_{2\nu+1}=\xi-(\nu+1)h,\ a_{2\nu+2}=\xi-(\nu+2)h,\ \cdots,\ a_n=\xi-(n-\nu)h.
\end{aligned}
\end{equation}
It requires $n\geq 2\nu$ and $a_n=\xi-(n-\nu)h>a$ to ensure the most information are used and all nodes are inside $[a,b]$.
This needs
\begin{eqnarray*}
	\frac{a-\xi+n(b-\xi)}{b-a}<\nu\leq\frac{n}{2}.
\end{eqnarray*}
Hence, once the number of nodes $n$ is conformed, the nodes $a_0,a_1,\cdots,a_n$ are conformed too.

Now, we can construct Lagrange interpolation polynomial at these nodes. Assume that $L_n(x)$ is the Lagrange polynomial which interpolates the function $f$ at $a_0,a_1,\cdots,a_n$ and written as
\begin{eqnarray*}
L_n(x)=\sum_{i=0}^nl_i(x)f(a_i),
\end{eqnarray*}
where $l_i(x)$ is the Lagrange fundamental polynomial mentioned in Theorem\ref{lem2.2}. The divided difference of
the interpolate polynomial $L_n(x)$ at the points $x$ and $\xi$ is obtained, where $\xi$ is repeated $p+1$ times,
\begin{eqnarray*}
L_n[x,\xi^{p+1}]=\sum_{i=0}^nl_i[x,\xi^{p+1}]f(a_i).
\end{eqnarray*}

Observe that $l_j(x)$ is a polynomial of degree $n$,  and the following property of the divided difference is an easy exercise:
\begin{equation*}
\begin{aligned}
l_i[x,\xi^{p+1}]&=\frac{1}{(x-\xi)^{p+1}}\Big(l_i(x)-\sum_{k=0}^p\frac{l^{(k)}_i(\xi)}{k!}(x-\xi)^k\Big)\\&=\frac{1}{(x-\xi)^{p+1}}\Big(l_i(\xi)+l'_i(\xi)(x-\xi)+\cdots+\frac{l_i^{(n)}(\xi)}{n!}(x-\xi)^n\\
&\ \ \ \ \ \ \ \ \ \ \ \ \ \ \  \ \ \ \ \ -\sum_{k=0}^p\frac{l^{(k)}_i(\xi)}{k!}(x-\xi)^k\Big)\\&
=\sum_{k=p+1}^n\frac{l^{(k)}_i(\xi)}{k!}(x-\xi)^{k-p-1}.
\end{aligned}
\end{equation*}
Therefore, the divided difference $f[x_c,\xi^{p+1}]$ can be rewritten as
\begin{eqnarray}\label{s3.3}
f[x_c,\xi^{p+1}]\approx  L_n[x_c,\xi^{p+1}]= \sum_{i=0}^nf(a_i)\sum_{k=p+1}^n\frac{l_i^{(k)}(\xi)}{k!}(x_c-\xi)^{k-p-1}.
\end{eqnarray}

\begin{theorem}\label{thm3.1}
Let $\displaystyle A_i^{(k)}=\frac{l_i^{(k)}(\xi)}{k!}$, nodes $a_0,a_1,\cdots,a_n$ are given by \eqref{s3.2},
then from Theorem \ref{lem2.2} we have
\begin{eqnarray}\label{s3.5}
A_{0}^{(k)}=Z_{k}(\eta_r),\ \ \ \eta_r=\frac{-1}{h^r}\Big[\sum_{j=1}^{\nu}\frac{1}{j^r}+\sum_{j=1}^{n-\nu}\frac{1}{(-j)^r}\Big].
\end{eqnarray}
For $i=1,2,\cdots,\nu$, the following formulae hold  
\begin{equation*}
\begin{aligned}
&A_{2i-1}^{(k)}=\frac{(-1)^{i-1}\nu!(n-\nu)!}{hi(\nu-i)!(n-\nu+i)!}Z_{k-1}\Big(\eta_r+\frac{1}{h^ri^r}\Big),\\
&A_{2i}^{(k)}=\frac{(-1)^{i}\nu!(n-\nu)!}{hi(\nu+i)!(n-\nu-i)!}Z_{k-1}\Big(\eta_r+\frac{1}{h^r(-i)^r}\Big),
\end{aligned}
\end{equation*}
and for $i=\nu+1,\nu+2,\cdots,n-\nu$, it follows as
\begin{equation*}
A_{\nu+i}^{(k)}=\frac{(-1)^{i}\nu!(n-\nu)!}{hi(\nu+i)!(n-\nu-i)!}Z_{k-1}\Big(\eta_r+\frac{1}{h^r(-i)^r}\Big).
\end{equation*}
\end{theorem}
\textbf{Proof}. Since $\xi=a_0$, this enables us to use equation \eqref{s2.5}, we have
\begin{eqnarray*}
A_0^{(k)}=\frac{l_0^{(k)}(\xi)}{k!}=Z_k(-S_{r0}(a_0)).
\end{eqnarray*}
From \eqref{s3.2}, the following equation can be easily verified
\begin{equation*}
\begin{aligned}
-S_{r0}(a_0)&=-\sum_{j=1}^n\frac{1}{(a_j-a_0)^r}\\
&=-\Big[\sum_{j=1}^{\nu}\Big(\frac{1}{(a_{2j-1}-a_0)^r}+\frac{1}{(a_{2j}-a_0)^r}\Big)+\sum_{j=2\nu+1}^{n}\frac{1}{(a_j-a_0)^r}\Big]\\
&=\frac{-1}{h^r}\Big[\sum_{i=1}^{\nu}\frac{1}{i^r}+\sum_{i=1}^{n-\nu}\frac{1}{(-i)^r}\Big]=\eta_r,
\end{aligned}
\end{equation*}
and this immediately infers \eqref{s3.5}.

For $i=1,2,\cdots,\nu$ and from \eqref{s2.7}$\Big(A_i^{(k)}=l'_i(\xi)Z_{k-1}(-S_{ri0}(\xi))\Big)$, it is easy to verify that
\begin{equation*}
\begin{aligned}
l_{2i-1}'(\xi)&=\frac{(-1)(-2)\cdots(-i+1)(-i-1)\cdots(-\nu)\cdot 1\cdot 2\cdots(n-\nu)h^{n-1}}{i(i-1)\cdots 1\cdot(-1)(-2)\cdots(i-\nu)(i+1)\cdots(i+n-\nu)h^{n}}\\
&=\frac{(-1)^{i-1}\nu!(n-\nu)!}{ih(\nu-i)!(n-\nu+i)!}.
\end{aligned}
\end{equation*}
Furthermore,
\begin{equation*}
\begin{aligned}
S_{r,2i-1,0}(\xi)&=\sum_{j=1\atop j\neq i}^{\nu}\frac{1}{(a_{2j-1}-a_0)^r}+\sum_{j=1}^{\nu}\frac{1}{(a_{2j}-a_0)^r}+\sum_{j=2\nu+1}^{n}\frac{1}{(a_j-a_0)^r}\\
&=\sum_{j=1\atop j\neq i}^{\nu}\frac{1}{(\xi+jh-\xi)^r}+\sum_{j=1}^{\nu}\frac{1}{(\xi-jh-\xi)^r}+\sum_{j=2\nu+1}^{n}\frac{1}{(\xi-(j-\nu)h-\xi)^r}\\
&=\frac{1}{h^r}\Big(\sum_{j=1\atop j\neq i}^{\nu}\frac{1}{j^r}+\sum_{j=1}^{n-\nu}\frac{1}{(-j)^r}\Big)=-\eta_r-\frac{1}{h^ri^r},
\end{aligned}
\end{equation*}
Together with $l'_{2i-1}(\xi)$, we get $A_{2i-1}^{(k)}$. In the case of $A_{2i}^{(k)}$ and $A_{\nu+i}^{(k)}$ can be obtained in a similar way. The proof is completed.\hfill $\square$

Combining  equation \eqref{s3.1}, \eqref{s3.3} and Theorem \ref{thm3.1}, we finally achieve the formula
\begin{eqnarray}\label{s3.5}
\begin{array}{ll}
\displaystyle H^*_{m,n,p}(\omega;f;\xi)=\sum_{k=1\atop k\neq
	c}^m\lambda_kf[x_k,\xi^{p+1}]+\lambda_c\sum_{i=0}^nf(a_i)\sum_{k=p+1}^nA_i^{(k)}(x_c-\xi)^{k-p-1}\\[.4cm]
\displaystyle \ \ \ \ \ \ \ \ \ \ \ \ \ \ \ \ \ \ \ \ \ \ \
+\sum_{j=0}^p\frac{f^{(j)}(\xi)}{j!}
\int_a^b\!\!\!\!\!\!\!\!\!=\frac{\omega(x)}{(x-\xi)^{p+1-j}}\mathrm{d}x.
\end{array}
\end{eqnarray}

It is worth to remark that the complexity of numerical quadrature proposed here is better than traditional method. It takes account of high order information and has the advantage that the coefficients could be calculated in a recurrent manner, which is usually stable and efficient in practical numerical calculations. Take $A_{2i-1}^{(k)}$ as example, we could calculate it in two steps. Firstly we calculate $$x_r:=\eta_r+\frac{1}{h^r i^r}=\frac{1}{h^r}\Big[\sum_{j=1}^{\nu}\frac{1}{j^r}+\sum_{j=1}^{n-\nu}\frac{1}{(-j)^r}+\frac{1}{i^r}\Big] (r=1,2,\cdots,k-1)$$ 
with $\displaystyle\sum_{r=1}^{k-1}(n+2)(r-1)=O(k^2 )$ multiplication operations(n is a fixed number). In the second step, the recurrence relation $\displaystyle(k-1)Z_{k-1}(x_r)=\sum_{j=k-1}^nx_jZ_{k-j-1}(x_r), Z_0=1$ with $1+2+\dots+k-1=O(k^2)$ multiplication operations is applied. Based on the above analysis, it is known that the solution of $A_{2i-1}^{(k)}$ requires only $O(k^2)$ operations. 

As the comparison, let us show that traditional scheme for calculating $\displaystyle A_i^{(k)}=\frac{l_i^{(k)}(\xi)}{k!}$ cost more than the one using cycle index. Notice that 
\begin{equation*}
	l_i^{(k)}(x)=\prod_{\substack{j=0 \\ j\neq i}}^{n} \frac{x-a_j}{a_i-a_j}=(\prod_{\substack{j=0 \\ j\neq i}}^{n} \frac{1}{a_i-a_j})\phi_n(x),\ \ \   \phi_n(x):=\prod_{\substack{j=0 \\ j\neq i}}^{n}(x-a_j).
\end{equation*}
According to Leibniz formula,
\begin{equation*}
\begin{aligned}
	\phi_n^{(k)}(x)&=\sum_{\substack{m_0+m_1+\dots+m_n=k \\ m_j\in \mathbb{N},\ \ j=0,1,\dots,i-1,i+1,\dots,n}}\frac{k!}{m_0!m_1!\dots m_n!}\prod_{\substack{j=0 \\ j\neq i}}^{n}(x-a_j)^{(m_j)}\\\\
	&=k!\sum_{0\leq j_0\leq j_1\leq\dots\leq j_{n-k}\leq n}\prod_{i=0}^{n-k}(x-a_{j_i}),
\end{aligned}
\end{equation*}
where we only consider the case that $m_j$ equals to 0 or 1(when $m_j\geq2$, $(x-a_j)^{(m_j)}\equiv 0$), thus the number of zero-valued $m_j$ is $(n-k+1)$, and their subscripts are denoted as $j_0,j_1,\dots,j_{n-k}$. Once all the combinations of $a_{j_i}$ are determined, we get the result of $A_i^{(k)}$. It's easy to verify that the above process requires $\binom{n+1}{k}(n-k)$ operations(n is a fixed number), which is much more than  $O(k^2)$.


\section{Error analysis}

In this section, we estimate the error of the quadrature formula in \eqref{s3.5}. The remainder term of the Gaussian quadrature formula \eqref{s3.1} with the nodes $x_1,x_2,\cdots,x_m$ is denoted as $R_{m} (g)$ where $g(x)=f[x,\xi^{p+1}]$, while the remainder term in polynomial interpolation with the nodes $a_0,a_1,\cdots,a_n$ is denoted as $R_{n} [x,\xi^{p+1}]$. The value $m$ and $n$ are given independently.

\begin{lemma}\label{lem4.1}

 Let $\Gamma$ be a simple closed curve in the complex plane surrounding the interval [-1,1] and D its interior. Having been subtracted out the singularity, the integrand g is an analytic function in D, thus the remainder term $R_{m} (g)$  admits the contour integral representation(cf. \cite{MC2003})
 \begin{eqnarray}\label{4.1}
 	\frac{1}{2\pi i}\oint_\Gamma K_m (z) g(z)\mathrm{d}z.
 \end{eqnarray}
 The kernel $K_m (z)$ is given by $K_m (z)=\dfrac{\rho_m (z)}{\pi_m (z)}$, where $\rho_m (z)=\int_{-1}^{1} \omega (t) \dfrac{\pi_m (t)}{z-t}\mathrm{d}t$, and $\pi_m (t)$ is the corresponding orthogonal polynomial with respect to the weight function $\omega (t)$ on $(-1,1)$.

\end{lemma}

\begin{lemma}\label{lem4.2}
Take the contour $\Gamma$ as a confocal ellipse with foci at the points -1, +1 and sum of semi-axes $\rho>1, \varepsilon_\rho=\{ z\in \mathbb{C}: z=\frac{1}{2} (\rho e^{i\vartheta} + \rho^{-1} e^{-i\vartheta} ), 0\leq\vartheta\leq2\pi \}$. D. B. Hunter(cf. \cite{BIT1995}) gives the expansion of $\rho_m (z)$ and $\pi_m (z)$ defined in \eqref{4.1}, together with an error estimation inequality as follows:
\begin{eqnarray}\label{4.2}
	\abs{R_{m} (g)} \leq \frac{4M(\varepsilon_\rho)\int_{-1}^{1} \omega (t)\mathrm{d}t}{\rho^{2m-1} (\rho-1)},
\end{eqnarray}
where $M(\varepsilon_\rho)=\max_{z\in \varepsilon_\rho} \abs{g(z)}$. Note that for certain specific weight functions, the above error bound can be further improved. For example, N. S. Kambo(cf. \cite{MC1970}) gives the error bound of Gauss-Legendre quarature formulae as $\abs{R_{m} (g)} \leq \dfrac{\pi M(\varepsilon_\rho) (\rho^2+1)}{\rho^{2m} (\rho^2-2)}$, which is better than \eqref{4.2} if $\rho$ is large ($\rho>\sqrt{2}$).
 
\end{lemma}

\begin{lemma}\label{lem4.3}
 Assume that $x_0,x_1,\cdots,x_k$ and $a_0,a_1,\cdots,a_n$ are sequences of nodes in $\mathbb{R}$ or in $\mathbb{C}$ ($k \leq n$). Then for any integer $m$ satisfying $0\leq m\leq k$, the divided difference of the remainder in polynomial interpolation can be expressed as (cf. \cite{JAT2004})
 \begin{small}
 \begin{eqnarray*}
 	R[x_0,x_1,\cdots,x_k]=\sum_{v=0}^{m-1} f[x_0,x_1,\cdots,x_v,a_0,a_1,\cdots,a_n] 
 	\omega_{n+1}[x_v,x_{v+1},\cdots,x_k] \notag\\ 
 	\shoveleft{+ \sum_{v=m}^{k} f[x_0,x_1,\cdots,x_v,a_0,a_1,\cdots,a_{n+m-v}]}{(x_v-a_{n+m-v}) \omega_{n+m-v}[x_v,x_{v+1},\cdots,x_k]},
 \end{eqnarray*} 
\end{small}
where $v=0,1,\cdots,n+1,$ and
 \begin{eqnarray*}
 	\omega_{v}(x):=\prod_{i=0}^{v-1}(x-a_i),\ \ \       R(x):=f(x)-\sum_{v=0}^{n}f[a_0,a_1,\cdots,a_{v}]\omega_{v}(x).
 \end{eqnarray*}

\end{lemma}

\begin{lemma}\label{lem4.4}
 The sequences of nodes $x_0,x_1,\cdots,x_k$ and $a_0,a_1,\cdots,a_n$ are defined in Lemma \ref{lem4.3}, then
 \begin{eqnarray*}
 	(f[x_0,x_1,\cdots,x_k,x])[a_0,a_1,\cdots,a_n]=f[x_0,x_1,\cdots,x_k,a_0,a_1,\cdots,a_n].
 \end{eqnarray*}
\textup{
We can use the definition of divided difference and apply traditional mathematics inductive approach to prove this lemma.}
 
\end{lemma}

\begin{theorem}\label{thm4.1}
Suppose that $f\in\mathbb{C}^{n+2}[-1,1]$, $f^{(n+1)}$,$f^{(n+2)}$ are bounded, and  $a_0,a_1,\cdots,a_n$ are given in \eqref{s3.2}, $n>p$. Meanwhile,  $ x_1,x_2,\cdots,x_m$ are chosen as the Gauss-Jacobi nodes with respect to the weight $\omega(t)=(1-t)^\alpha(1+t)^\beta$ on $(-1,1)$ and $\rho,M(\varepsilon_\rho)$ are defined in Lemma \autoref{lem4.2}. Then we have the error estimation formula
 
\begin{equation*}
	\abs{R_{m,n}^{*}}\leq \frac{2^{\alpha+\beta+3} \Gamma(\alpha+1) \Gamma (\beta+1) M(\varepsilon_\rho)}{\Gamma(\alpha+\beta+2) \rho^{2m-1} (\rho-1)}+\frac{M_1+ M_2}{p!}\frac{1}{(n-p)!}(1+\frac{p+2}{n})^{n-p}.
\end{equation*}
 where $M_1=\max\abs{f^{(n+1)}}, M_2=\max\abs{f^{(n+2)}}$.\\
\end{theorem} 
\textup{
\textbf{Proof}. Substituting the weight function $\omega(t)=(1-t)^\alpha(1+t)^\beta$ into \eqref{4.2}, we obtain that}
\begin{equation*}
	\abs{R_{m} (g)}\leq \frac{2^{\alpha+\beta+3} \Gamma(\alpha+1) \Gamma (\beta+1) M(\varepsilon_\rho)}{\Gamma(\alpha+\beta+2) \rho^{2m-1} (\rho-1)},
\end{equation*}
In practice, $M(\varepsilon_\rho)$ is calculated for several values of $\rho$ within the appropriate range and the smallest number so obtained is taken as the upper bound for $R_{m} (g)$. Once $\alpha,\beta,\rho$ are fixed and $m\to$ $\infty$, we get $R_{m} (g)\to 0$.

According to the definition of the remainder term in polynomial interpolation, 
\begin{equation*}
 	R_{n} [x,\xi^{p+1}]=f[x,\xi^{p+1}]-\sum_{v=0}^{n}f[a_0,a_1,\cdots,a_v]\omega_{v}[x,\xi^{p+1}].
\end{equation*}
From Lemma \ref{lem4.3},  $\alpha_0,\beta_0,\alpha_v\in (-1,1)$, we know that
 
\begin{equation*}
 	\begin{split}
 		R_{n} [x,\xi^{p+1}]&=R[x,\xi^{p+1}]= f[x,a_0,a_1,\cdots,a_n] 
 		\omega_{n+1}[x,\xi^{p+1}] \\
 		&+ \sum_{v=1}^{p+1} f[x,\xi^{v},a_0,a_1,\cdots,a_{n+1-v}]
 		(\xi-a_{n+1-v}) \omega_{n+1-v}[\xi^{p+2-v}]\\
 		&=\frac{f^{(n+1)}(\alpha_0)}{(n+1)!}\frac{\omega_{n+1}^{(p+1)}(\beta_0)}{(p+1)!}+\sum_{v=1}^{p+1}\frac{f^{(n+2)}(\alpha_v)}{(n+2)!} (\xi-a_{n+1-v})\frac{\omega_{n+1-v}^{(p+1-v)}(\xi)}{(p+1-v)!}.
 	\end{split} 
\end{equation*}
\textup{
From Leibniz formula, 
\begin{equation*}
\begin{aligned}
 	\omega_{n+1}^{(k)}(x)&=\sum_{\substack{m_0+m_1+\dots+m_n=k \\ m_j\in \mathbb{N},\ \ j=0,1,\dots,n}}\frac{k!}{m_0!m_1!\dots m_n!}\prod_{j=0}^{n}(x-a_j)^{(m_j)}\\
 	&=k!\sum_{0\leq j_0\leq j_1\leq\cdots\leq j_{n-k}\leq n }\prod_{i=0}^{n-k}(x-a_{j_i}).
\end{aligned}
\end{equation*}
Let $a_{j_i}=\xi-\lambda_{j_i}h, h=\frac{1-\xi}{\nu+1}$ (defined in \eqref{s3.2}), where} \\
\begin{equation*}
 	 -\nu\leq\lambda_{j_i}\leq n-\nu,\\\\
 	 \nu>\frac{-1-\xi+n(1-\xi)}{2}, 
\end{equation*}
\textup{
thus we obtain $h < \frac{2}{n+1}$, and then}
\begin{equation*}
	\begin{split}
		\prod_{i=0}^{n-k}(x-a_{j_i})
		&\leq \Big(\frac{\sum_{i=0}^{n-k}(x-\xi+\lambda_{j_i}h)}{n-k+1}\Big)^{n-k+1}\\
		&=\Big(\frac{(n-k+1)(x-\xi)+(\sum_{i=0}^{n-k}\lambda_{j_i})h}{n-k+1}\Big)^{n-k+1}\\
		&\leq \Big(\frac{(n-k+1)(x-\xi)+(\sum_{\lambda_{j_i}=k-\nu}^{n-\nu}\lambda_{j_i})h}{n-k+1}\Big)^{n-k+1}\\
 		&=\Big(\frac{(n-k+1)(x-\xi)+\frac{(n+k-2\nu)(n-k+1)}{2}h}{n-k+1}\Big)^{n-k+1}\\
 		&\leq\Big (1-\xi+\frac{n+k-2\nu}{2}h\Big)^{n-k+1}   
 		=\Big(\nu+1+\frac{(n+k-2\nu)}{2}\Big)^{n-k+1}h^{n-k+1}\\
 		&=\Big(\frac{n+k+2}{2}\Big)^{n-k+1} h^{n-k+1}\leq \Big(1+\frac{k+1}{n+1}\Big)^{n-k+1},
 	\end{split} 	
\end{equation*}
\textup{
so we get $\omega_{n+1}^{(k)}(x)\leq k!\binom{n+1}{k} (1+\frac{k+1}{n+1})^{n-k+1}$.
Therefore,}
\begin{equation*}
\begin{aligned}
 		& \abs{R_{n}[x,\xi^{p+1}]} \\
 		\hspace{3mm}& \leq\frac{\abs{f^{(n+1)}(\alpha_0)}}{(n+1)!}\frac{\omega_{n+1}^{(p+1)}(\beta_0)}{(p+1)!}+\sum_{v=1}^{p+1}\frac{\abs{f^{(n+2)}(\alpha_v)}}{(n+2)!} (\xi-a_{n+1-v})\frac{\omega_{n+1-v}^{(p+1-v)}(\xi)}{(p+1-v)!} \\
 		& \leq \frac{M_1}{(n+1)!}\frac{(n+1)!}{(p+1)!(n-p)!}(1+\frac{p+2}{n+1})^{n-p}  
 		+\sum_{v=1}^{p+1}\frac{M_2}{(n+2)!} (n-\nu)h \binom{n+1-v}{p+1-v}(1+\frac{p+2-v}{n+1-v})^{n-p}\\
 		& \leq \frac{M_1}{(n+1)!}\frac{(n+1)!}{(p+1)!(n-p)!}(1+\frac{p+2}{n+1})^{n-p} 
%
		+ \sum_{v=1}^{p+1}\frac{M_2}{(n+2)!} (n-\nu)h \binom{n+1-v}{p+1-v}(1+\frac{p+2}{n+1})^{n-p}\\
		&\leq
		\frac{M_1}{(n+1)!}\frac{(n+1)!}{(p+1)!(n-p)!}(1+\frac{p+2}{n+1})^{n-p} 
		+\frac{M_2}{(n+2)!} (n-\nu) \frac{2}{n+1} \binom{n+1}{p}(1+\frac{p+2}{n+1})^{n-p}\\
 		&\leq \frac{M_1}{p!}\frac{1}{(n-p)!}(1+\frac{p+2}{n+1})^{n-p}+\frac{ M_2}{p!} \frac{1}{(n-p)!}(1+\frac{p+2}{n+1})^{n-p}\\
 		&=\frac{M_1+ M_2}{p!}\frac{1}{(n-p)!}(1+\frac{p+2}{n+1})^{n-p}.  
\end{aligned}
\end{equation*}
\textup{
Once p is fixed and n$\to$ $\infty$, $(1+\frac{p+2}{n+1})^{n-p}\to$ ${\textit{e}}^{p+2}$, where $\textit{e}$ is the base of natural logarithms, thus we get $R_{n}[x,\xi^{p+1}]\to 0$. Therefore, the convergence of the error estimation formula has been proved. \hfill $\square$}

\section{Numerical {\color{huxl}Examples}}

 In this section, we present some numerical examples to demonstrate the applicability and accuracy of the quadrature rule proposed in this
 paper. In the quadrature formula \eqref{s3.5}, $ x_1,x_2,\cdots,x_m$ are chosen as the Gaussian nodes with respect to the weight $\omega(x)$. Furthermore, the formula involves another group of nodes $a_0,a_1,\cdots,a_n$. Here absolute error is defined as the absolute value of the true error and the missing points in the figures indicate that the data have reached machine precision. Note that all computations carried out in this paper were performed using the package of MATLAB 2017b on a
 personal computer with the operating system Windows 10, 8 GB RAM and 2.11 GHz CPU clock speed.
 
\begin{example}\label{ex1}
Here, we give the following Hadamard finite-part integral to verify some properties of the new quadrature rule.
\begin{eqnarray*}
I_1(\xi;p)=\int_{-1}^1\!\!\!\!\!\!\!\!\!\!=\
\frac{e^x}{(x-\xi)^{p+1}}\mathrm{d}x,
\end{eqnarray*}
where $x_1,x_2,\cdots,x_m$ are zeros of Legendre polynomial.
\end{example}

\begin{figure}[h]
	\centering
	\includegraphics[width=10cm,height=6cm]{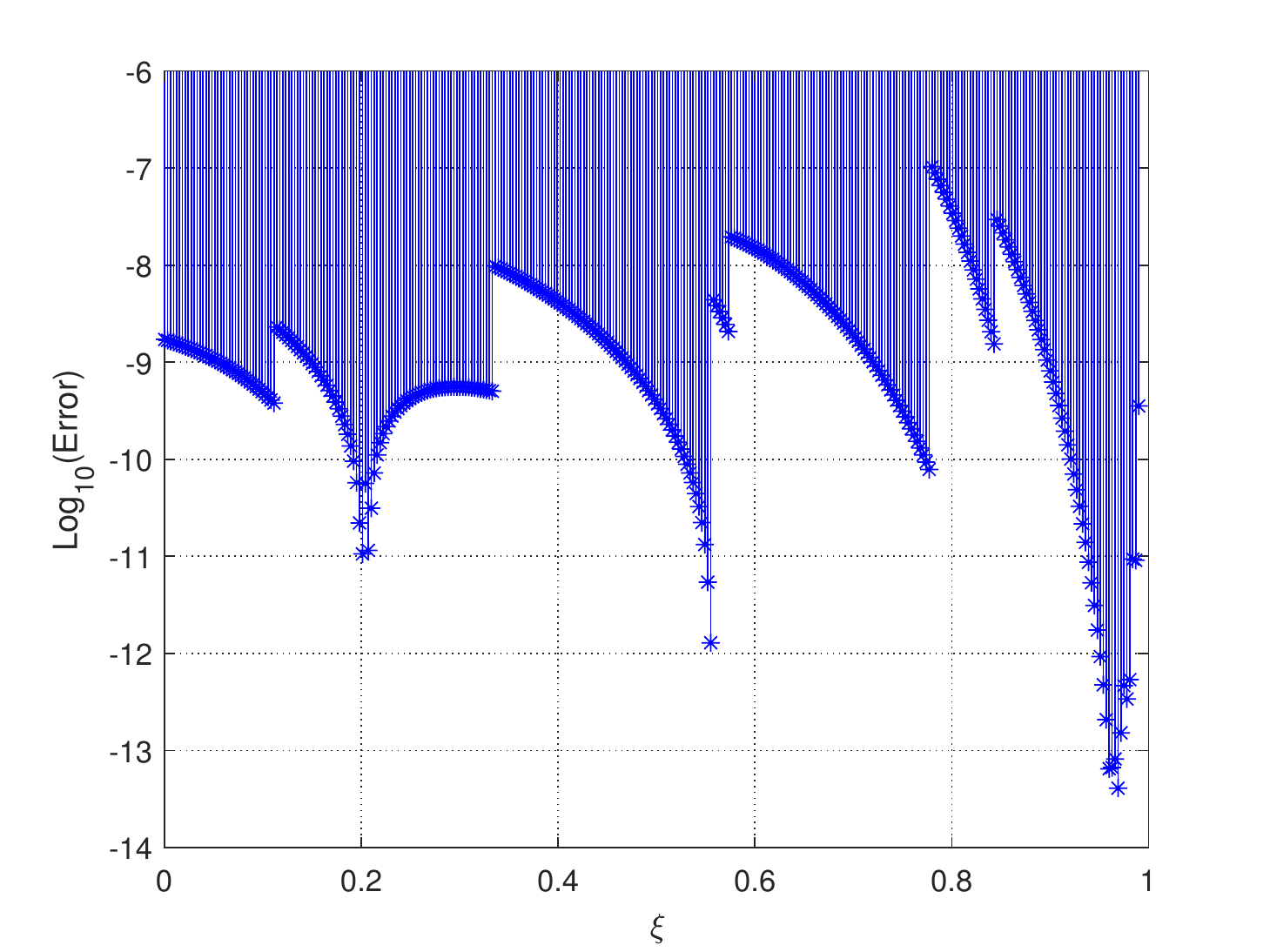}
	\caption{the absolute errors of the approximations to $I_1(p=0,m=7,n=8)$ for different singularities $\xi$} 
	\label{fig:1} 
\end{figure}

If $p=0$, the exact solution of the above integral can be expressed as
\begin{equation*}
	[Ei(1-\xi)-Ei(-1-\xi)]{\textit{e}}^{ \xi},
\end{equation*}
where $Ei(x)=\int_{-\infty}^{x} \frac{{\textit{e}}^t}{t} \mathrm{d}t$ represents the exponential integral. Once we fix the value of $m$ and $n$, the absolute errors of the approximations to $I_{1}(\xi;p)$ can be obtained for different singularities. From \autoref{fig:1}, we observe that the error is bounded uniformly, namely, independently of $\xi$, which is an illustrative evidence of uniform error bound mentioned in Section 4.

\begin{figure}[h]
	\centering
	\subfigure
	{
		\begin{minipage}[b]{.47\linewidth}
			\centering
			\includegraphics[width=7.48cm,height=5cm]{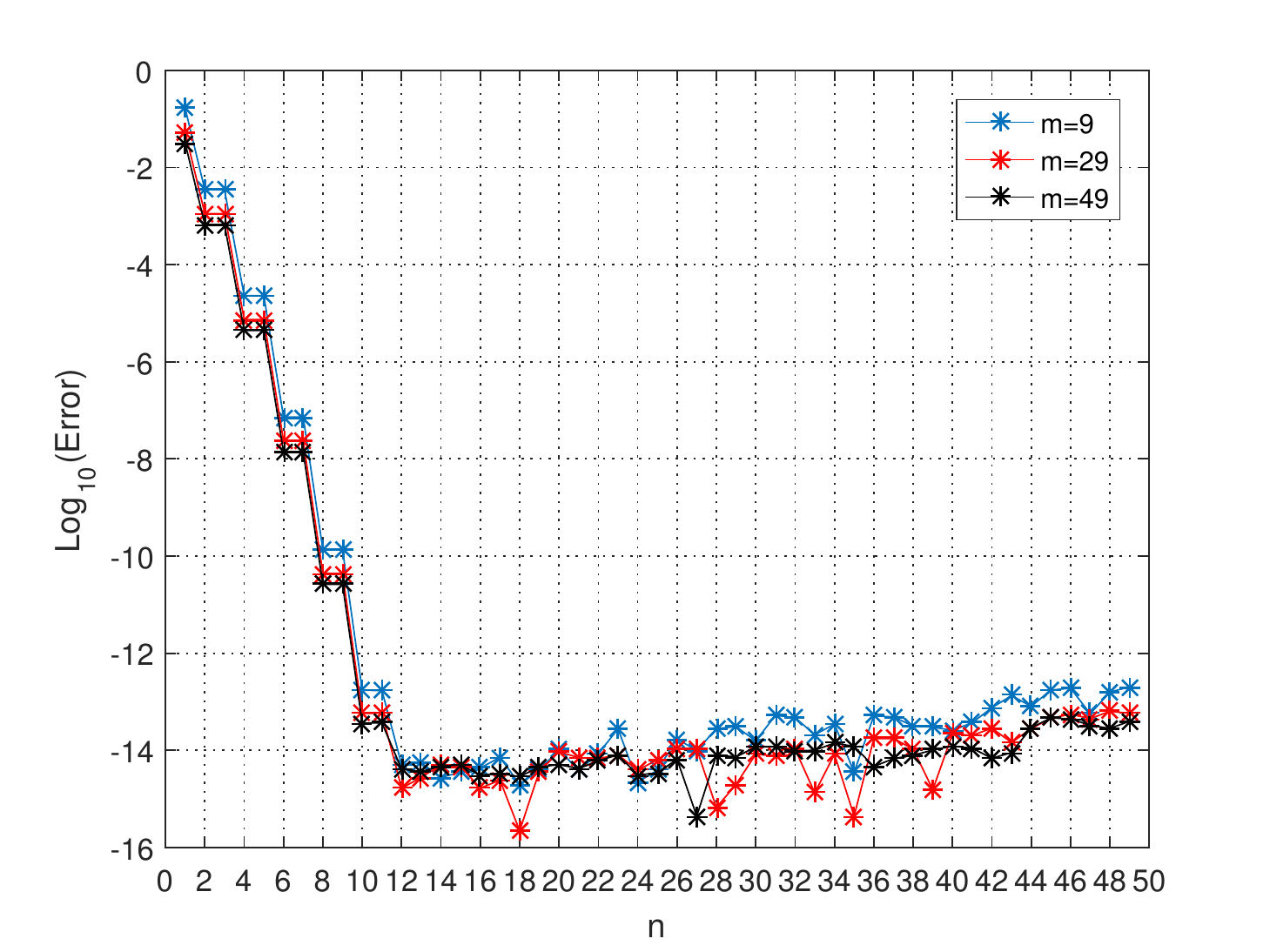}
		\end{minipage}
	}
	\subfigure
	{
		\begin{minipage}[b]{.47\linewidth}
			\centering
			\includegraphics[width=7.48cm,height=5cm]{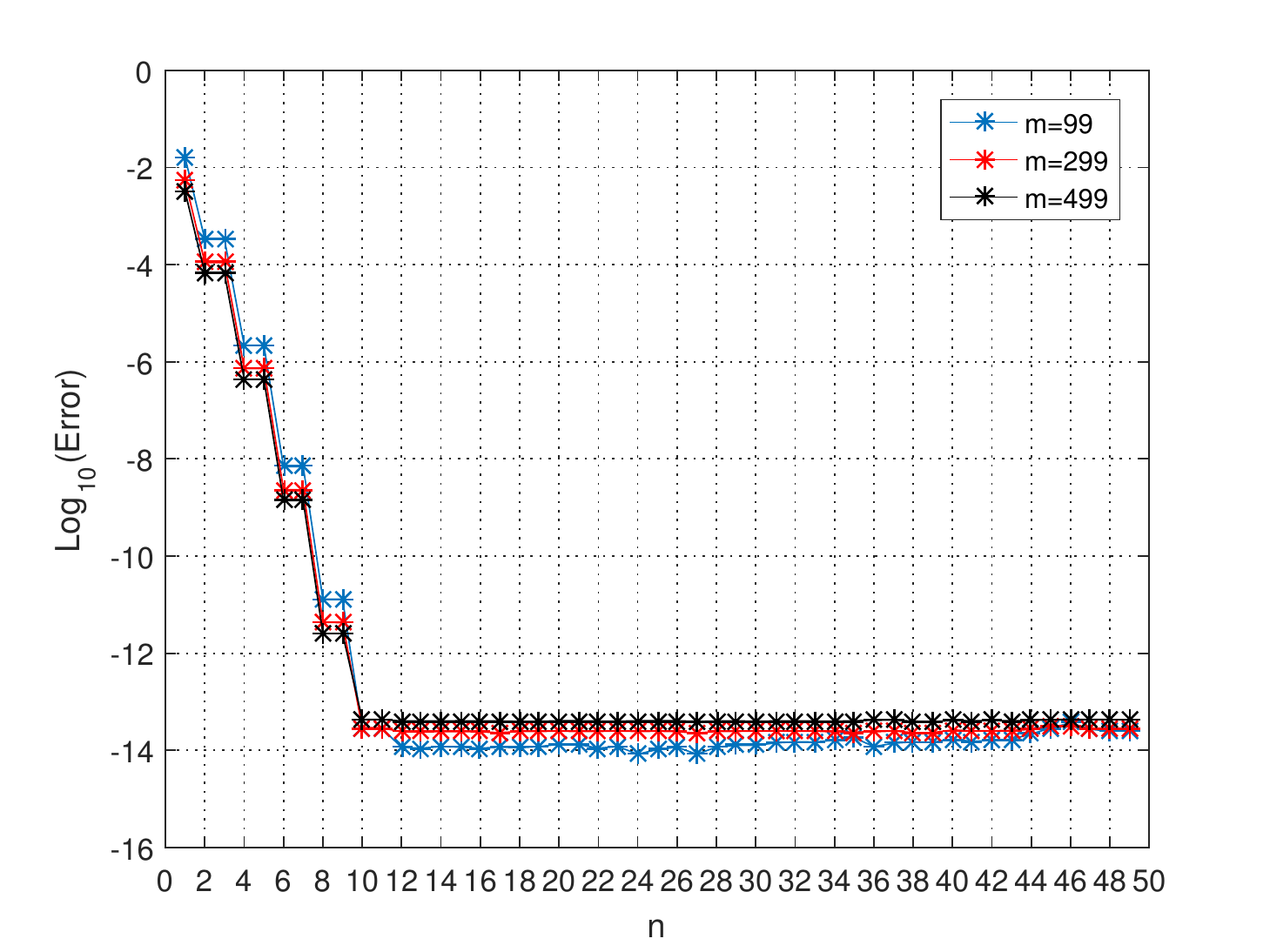}
		\end{minipage}
	}
	\caption{the absolute errors of the approximations to $I_1(p=1)$ for diffrent $n$ and  $m$ }
	\label{fig:2}
\end{figure}

\begin{figure}[h]
	\centering
	\subfigure
	{
		\begin{minipage}[b]{.47\linewidth}
			\centering
			\includegraphics[width=7.48cm,height=5cm]{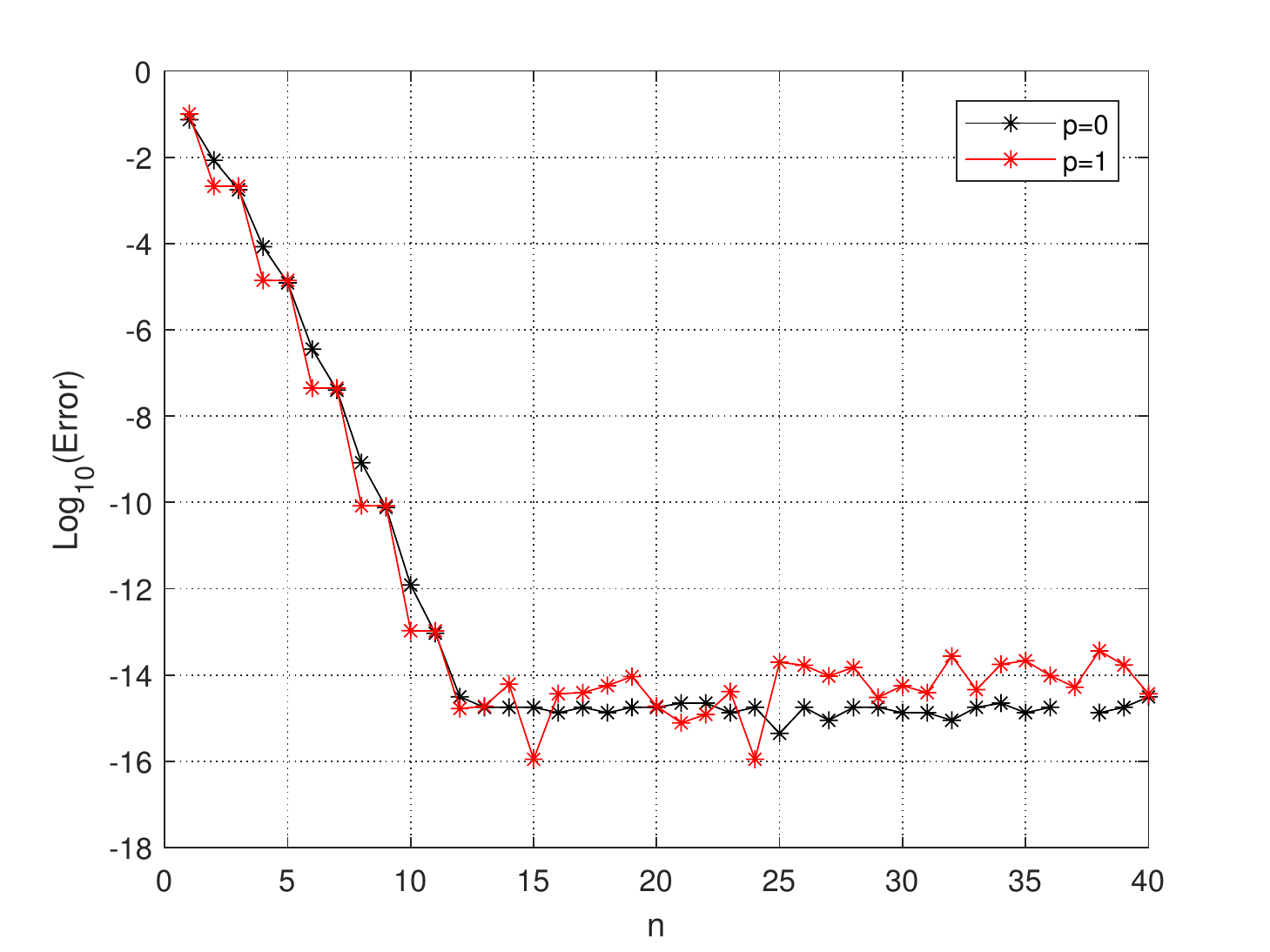}
		\end{minipage}
	}
	\subfigure
	{
		\begin{minipage}[b]{.47\linewidth}
			\centering
			\includegraphics[width=7.48cm,height=5cm]{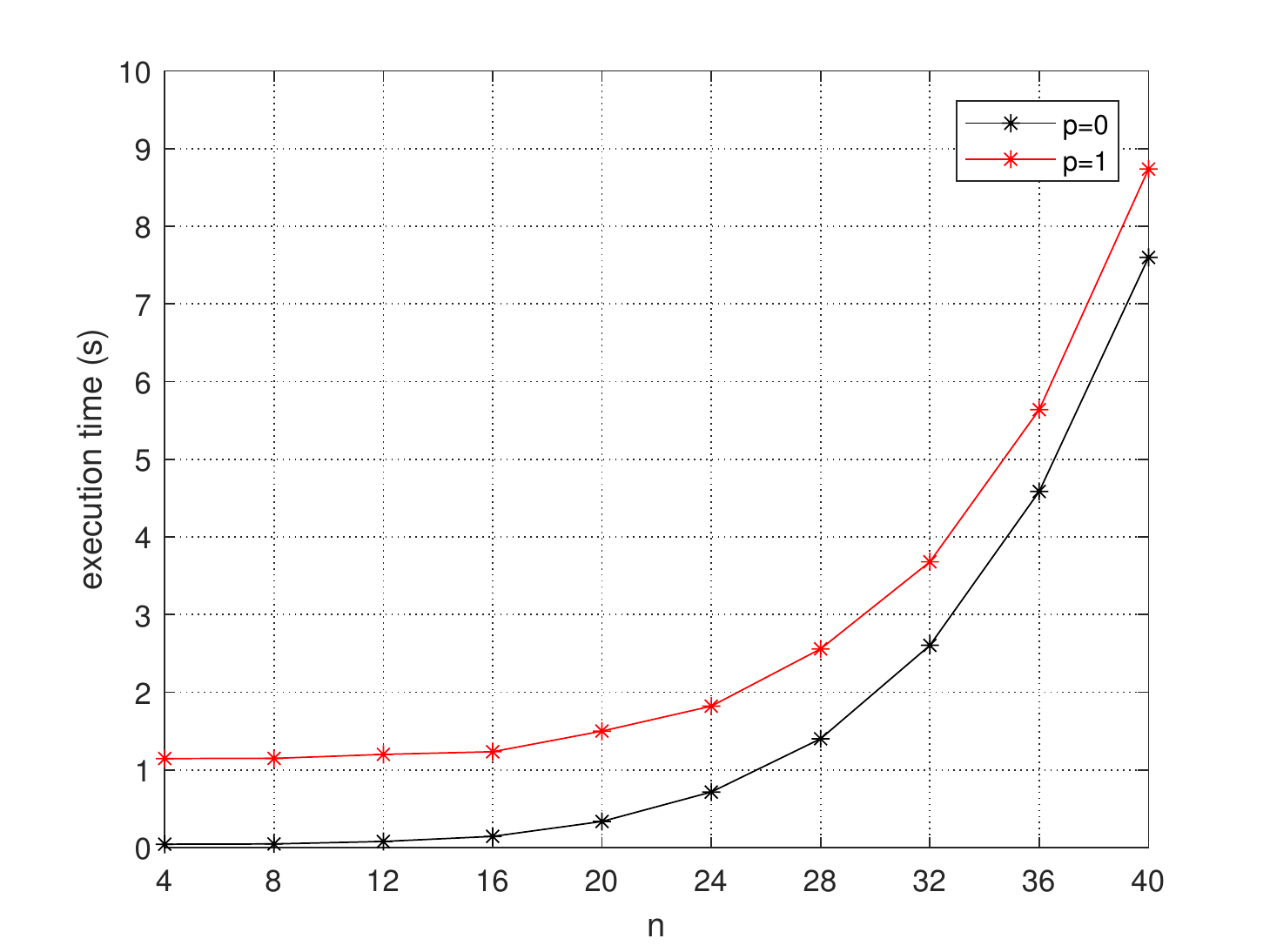}
		\end{minipage}
	}
	\caption{the absolute errors and execution time of the approximations to $I_1(m=15, \xi=10^{-5})$ for different value $p$}
	\label{fig:3}
\end{figure}

When $\xi=10^{-5}$, we present the absolute errors $R_{m,n}^*$ of $I_1$ in \autoref{tab1} with varied values of $n$ and different numbers of
quadrature nodes $m$, where $p$ is given by $0$ and $1$. Meanwhile, we also record the execution time of our algorithm for each m and n (take the average of ten measurements).
From this table, we can see that for both $p=0$ and $p=1$, the error decreases sharply at the beginning while $n$ is increasing but when $n$ is large enough the error fluctuates up and down in a small range or even rises up. We can also see this phenomenon on the left figure of \autoref{fig:2}. Meanwhile, a bigger $n$ will lead to longer execution time for our algorithm(see \autoref{fig:3}). 
\begin{table}[h]
	\small\centering
	\caption{the absolute errors $R_{m,n}^*$ of $I_1(\xi=10^{-5})$ for $p=0$ and $p=1$}\label{tab1}
	\setlength{\tabcolsep}{3mm}
		\begin{tabular}{c c c| c c }
			\hline
			& $p=0$ & & $p=1$  &  \\ \hline
			$\mathrm{n}$& $R_{m,n}^*(m=7)$ & $R_{m,n}^*(m=15)$ & $R_{m,n}^*(m=7)$ & $R_{m,n}^*(m=15)$   \\
			\hline
			4 & $1.742840 \mathrm{ E}-04$ & $8.447273 \mathrm{ E}-05$ & $2.897117 \mathrm{ E}-05$&$1.404187 \mathrm{ E}-05$ \\
			8 & $1.716980 \mathrm{ E}-09$ & $8.321455 \mathrm{ E}-10$& $1.719901 \mathrm{ E}-10$& $8.335765 \mathrm{ E}-11$\\
			11 & $1.940670 \mathrm{ E}-13$ & $9.281464 \mathrm{ E}-14$& $2.134959 \mathrm{ E}-13$& $1.051381 \mathrm{ E}-13$\\
			12 & $3.108624 \mathrm{ E}-15$ & $3.108624 \mathrm{ E}-15$& $6.328271 \mathrm{ E}-15$ & $1.665335 \mathrm{ E}-15$\\
			24 & $4.440892 \mathrm{ E}-16$ & $1.776357 \mathrm{ E}-15$& $1.765254 \mathrm{ E}-14$& $1.110223 \mathrm{ E}-16$\\
			31 & $2.664535 \mathrm{ E}-15$ &$1.332267 \mathrm{ E}-15$& $4.107825 \mathrm{ E}-15$& $3.885781 \mathrm{ E}-15$\\
			44 & $2.220446 \mathrm{ E}-15$ &$2.220446 \mathrm{ E}-15$& $1.221245 \mathrm{ E}-15$ & $2.409184 \mathrm{ E}-14$\\
			59 & $3.996803 \mathrm{ E}-15$ &$1.776357 \mathrm{ E}-15$& $1.282308 \mathrm{ E}-13$ & $6.694645 \mathrm{ E}-14$\\
			\hline
		\end{tabular}
\end{table}

 Therefore we're supposed to find an optimal option of $n=\hat{n}$ which makes the error relatively small. In  \autoref{fig:2}, we find that when taking different values of $m$,  the value of $n$ tends to be stable.  In particular when $m$ takes a larger number, the absolute errors  increase. To reduce the execution time effectively, we can take the above stable value $n$ as an optimal value $\hat{n}$ , especially when $m$ is relatively large. For small values of $m$, like those less than or closer to the stable value $n$, we can limit the search to a certain small range $[m-\delta,m+\delta]$.
For instance, when $m=15,\delta=10$ for $p=0$ and $p=1$ it's easy to calculate by MATLAB that $\hat{n}$ equals to 25 and 15, respectively (see \autoref{fig:3}). 

\begin{example}\label{ex2}
In this example, we show numerical results for the following Hadamard finite-part integral:
\begin{eqnarray*}
I_2=\int_{-1}^1\!\!\!\!\!\!\!\!\!\!=\
\frac{(1.21-x^2)^{-1/2}}{(x-10^{-5})^{2}}\mathrm{d}x.\ \ \
\end{eqnarray*}
In this case $\omega(x)=1$, so that $x_m$ are zeros of Legendre polynomial and the exact solution of $I_2$ is $-0.757450528292818$.
\end{example}

Firstly we make certain an appropriate value $n=\hat{n}$ for each m. It can be discovered from \autoref{fig:4} that for a fixed $m$, the error curve will not decline continuously when $n$ is large enough, thus $\hat{n}$ can be determined using the method mentioned in the last experiment.
\begin{figure}[h]
	\centering
	\subfigure
	{
		\begin{minipage}[b]{.47\linewidth}
			\centering
			\includegraphics[width=7.48cm,height=5cm]{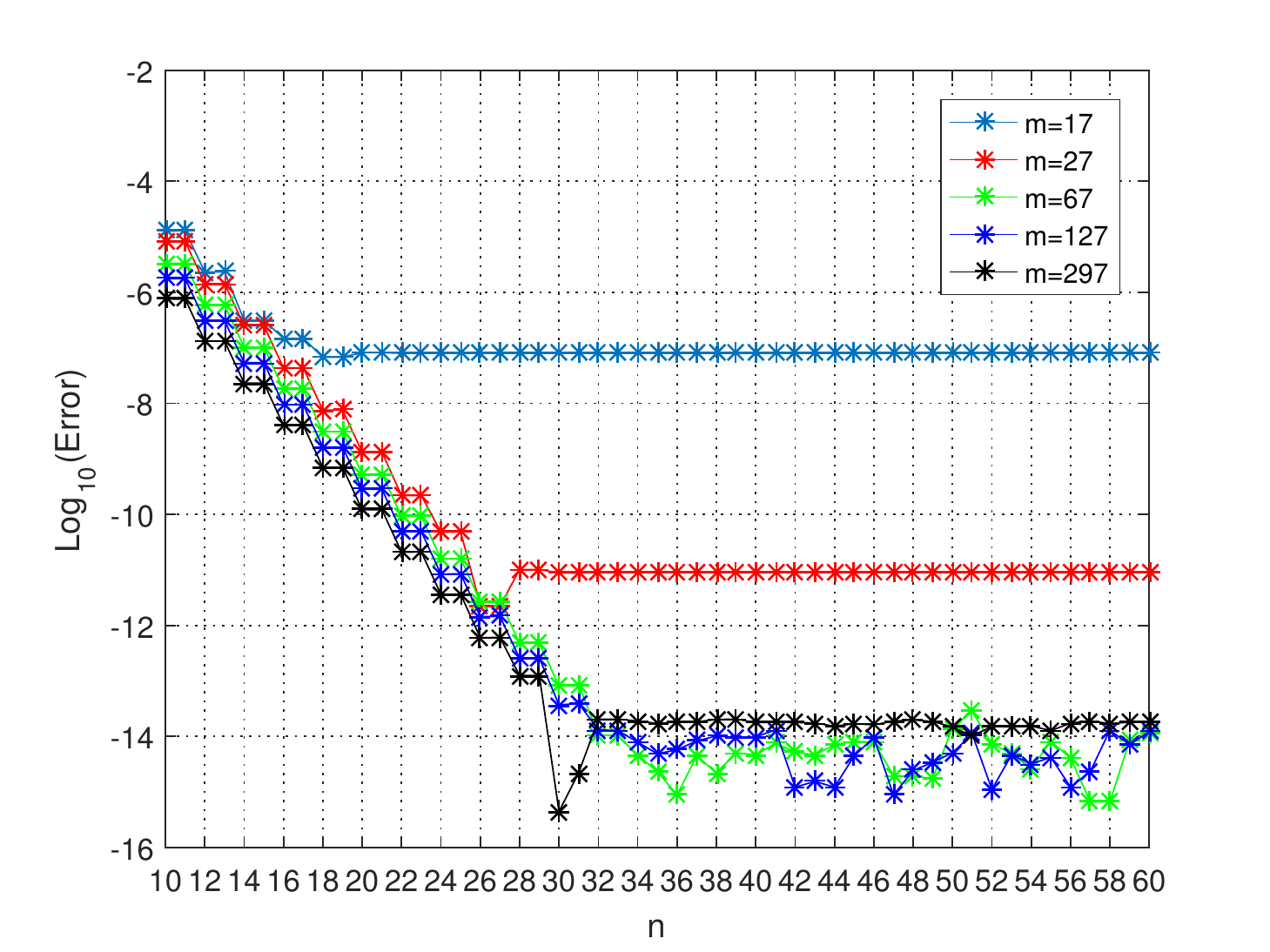}
		\end{minipage}
	}
	\subfigure
	{
		\begin{minipage}[b]{.47\linewidth}
			\centering
			\includegraphics[width=7.48cm,height=5cm]{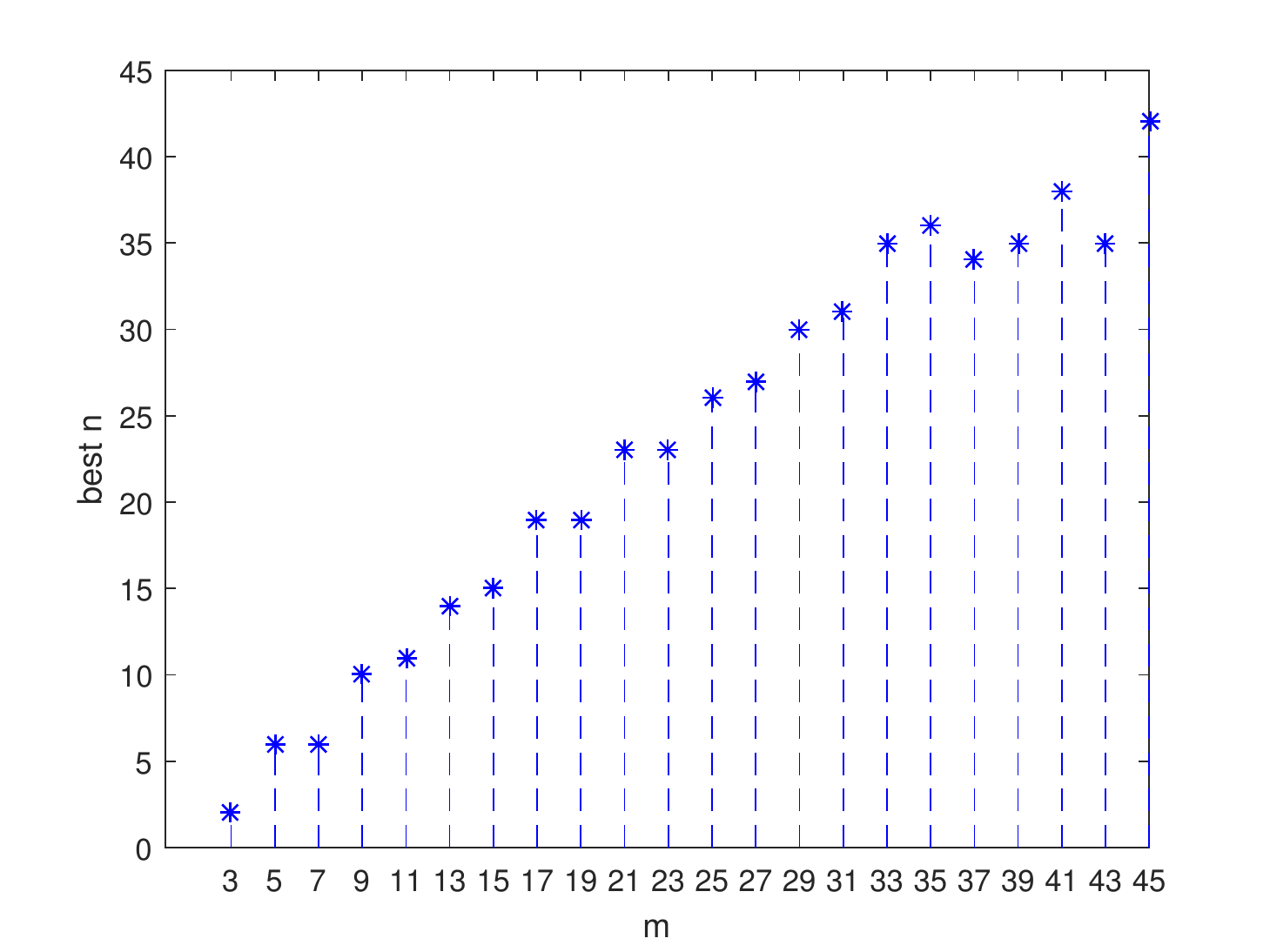}
		\end{minipage}
	}
	\caption{the absolute errors of the approximations to $I_2$ for diffrent $n$ and the best option of $n=\hat{n}$ for each $m$}
	\label{fig:4}
\end{figure}

Then we compare numerically the quadrature formula(\ref{s3.5}) (denoted as Algorithm 3) with the Gaussian quadrature formula considered in \cite{MC85} and another quadrature formula proposed in \cite{JCAM97}(denoted as Algorithm 1 and Algorithm 2 respectively). Some results are shown in \autoref{tab2}. We see that our new quadrature rule converges much faster and has higher rate of convergence. The only disadvantage is the longer execution time. The main reason for the longer time is that it takes more time to search the optimal $n$ value. But we can reduce it efficiently by using the stable value $n$ for bigger $m$, which might not be the best option but makes the absolute error have the same order of magnitude as the best situation. 

\begin{table}[ht]
	\small
	\centering
\caption{Comparison of our method with Gaussian quadrature formula considered in \cite{MC85} and
another quadrature formula proposed in \cite{JCAM97} for integrals $I_2$}\label{tab2}
    \begin{tabular}{c c c c c c c}
    \hline
 & Algorithm 1 & & Algorithm 2  & & Algorithm 3 &  \\ \hline
$\mathrm{m}$& absolute error & time(s) & absolute error &time(s) & absolute error & time(s)   \\
\hline
 3 & $3.926491 \mathrm{ E}-02$ & $0.4991$ & $1.768617 \mathrm{ E}-01$ & $0.0851$ & $2.336930 \mathrm{ E}-02$ & $0.2786$  \\
 9 & $1.239773 \mathrm{ E}-04$ & $0.6375$& $9.282299 \mathrm{ E}-03$& $0.0862$ & $1.047047 \mathrm{ E}-04$& $0.7463$\\
 15 & $1.541212 \mathrm{ E}-07$ & $0.7231$& $4.928689 \mathrm{ E}-04$& $0.0868$ & $5.377683 \mathrm{ E}-08$& $1.2668$\\
 21 & $2.980971 \mathrm{ E}-06$ &$0.8336$& $2.900791 \mathrm{ E}-05$& $0.0870$ &$1.792015 \mathrm{ E}-09$& $2.1794$\\
 27 & $4.032681 \mathrm{ E}-06$ & $1.0016$& $1.781499 \mathrm{ E}-06$& $0.1022$ & $2.230882 \mathrm{ E}-12$& $2.9998$\\
 33 & $3.435934 \mathrm{ E}-06$ &$1.2002$& $1.123160 \mathrm{ E}-07$& $0.1215$ &$4.485301 \mathrm{ E}-14$& $6.8992$\\
 39 & $1.207606 \mathrm{ E}-06$ &$1.4634$& $7.204918 \mathrm{ E}-09$& $0.1331$ &$9.325873 \mathrm{ E}-15$& $12.3193$\\
 45 & $3.254382 \mathrm{ E}-06$ &$2.4647$& $4.679410 \mathrm{ E}-10$& $0.1484$ &$7.105427 \mathrm{ E}-15$& $18.1009$\\
  \hline
    \end{tabular}
\end{table}

\begin{example}\label{ex3}
In this example, we show numerical results for the following Hadamard finite-part integral:
\begin{eqnarray*}
I_3(\xi,\lambda)=\int_{-1}^1\!\!\!\!\!\!\!\!\!\!=\
\frac{(x^2+\lambda^2)^{-1}}{(x-\xi)^2}\frac{\mathrm{d}x}{\sqrt{1-x^2}},\ \ \ \ \xi=0.25.
\end{eqnarray*}
In this case $\omega(x)=1/\sqrt{1-x^2}$, so that $x_m$ are zeros of Chebyshev polynomial of the first kind and the exact solution of $I_3(\xi,\lambda)$ is
\begin{eqnarray*}
	\frac{\pi (\xi^2-\lambda^2)}{\lambda \sqrt{\lambda^2+1} (\lambda^2+\xi^2)^2}.
\end{eqnarray*}
\end{example}

In this example, we present the absolute errors $R_m^*$ of $I_3$ in \autoref{tab3} with varied values of $\lambda$ (so that the integrand functions differ from each other). It can be discovered that for a bigger value $\lambda$ we have a higher rate of convergence in our quadrature formula. Then we compare numerically our method with the quadrature formula proposed in \cite{JCAM97} and show the results in \autoref{fig:4}. For varied values of $\lambda$, the quadrature formula \eqref{s3.5} is of higher accuracy than the other method.
\begin{table}[h]
	\centering
	\small
	\caption{the absolute errors $R_m^*$ of $I_4(\xi=0.25)$ for different values of $\lambda$ }\label{tab3}
	\setlength{\tabcolsep}{1mm}
		\begin{tabular}{c c c c c c c} 
			\hline
			& $\lambda=1.5$ & & $\lambda=2.5$ & &$\lambda=5.0$&
			\\
			\hline
			$\mathrm{m}\ \ \ $ & $\hat{n}$ & $R^*_m$ & $\hat{n}$ & $R^*_m$ & $\hat{n}$ & $R^*_m$ \\
			\hline
			3 \ \ \ & 5& $4.638894\mathrm{ E}-04$ &10& $7.395720\mathrm{ E}-06$ &10& $9.222675\mathrm{ E}-09$\\
			4 \ \ \ & 6& $2.586607\mathrm{ E}-06$ & 6& $2.677880\mathrm{ E}-07$ &10& $9.052034\mathrm{ E}-11$\\
			5 \ \ \ &14& $6.152098\mathrm{ E}-06$ &19& $9.931182\mathrm{ E}-09$ &16& $1.001525\mathrm{ E}-13$\\
			6 \ \ \ &10& $2.894122\mathrm{ E}-07$ &10& $1.404692\mathrm{ E}-11$ &12& $8.596422\mathrm{ E}-15$\\
			7 \ \ \ &12& $2.906012\mathrm{ E}-08$ &13& $4.271375\mathrm{ E}-12$ &14& $4.796510\mathrm{ E}-16$\\
			8 \ \ \ &20& $4.742672\mathrm{ E}-09$ &18& $5.149353\mathrm{ E}-13$ &23& $4.510281\mathrm{ E}-17$\\
			9 \ \ \ &21& $4.349271\mathrm{ E}-10$ &21& $9.450773\mathrm{ E}-15$ &15& $1.734723\mathrm{ E}-17$\\
			10 \ \ \ &21& $4.008799\mathrm{ E}-11$ &21& $4.857226\mathrm{ E}-16$ &15& $4.423545\mathrm{ E}-17$\\
			11 \ \ \ &22& $3.663181\mathrm{ E}-12$ &22& $6.661338\mathrm{ E}-16$ &12& $1.908196\mathrm{ E}-17$\\
			12 \ \ \ &21& $4.963807\mathrm{ E}-13$ &20& $9.992007\mathrm{ E}-16$ &17& $8.073617\mathrm{ E}-19$\\
			\hline
		\end{tabular}
\end{table}

\begin{figure}[h]  
	\centering
	\includegraphics[width=10cm,height=7cm]{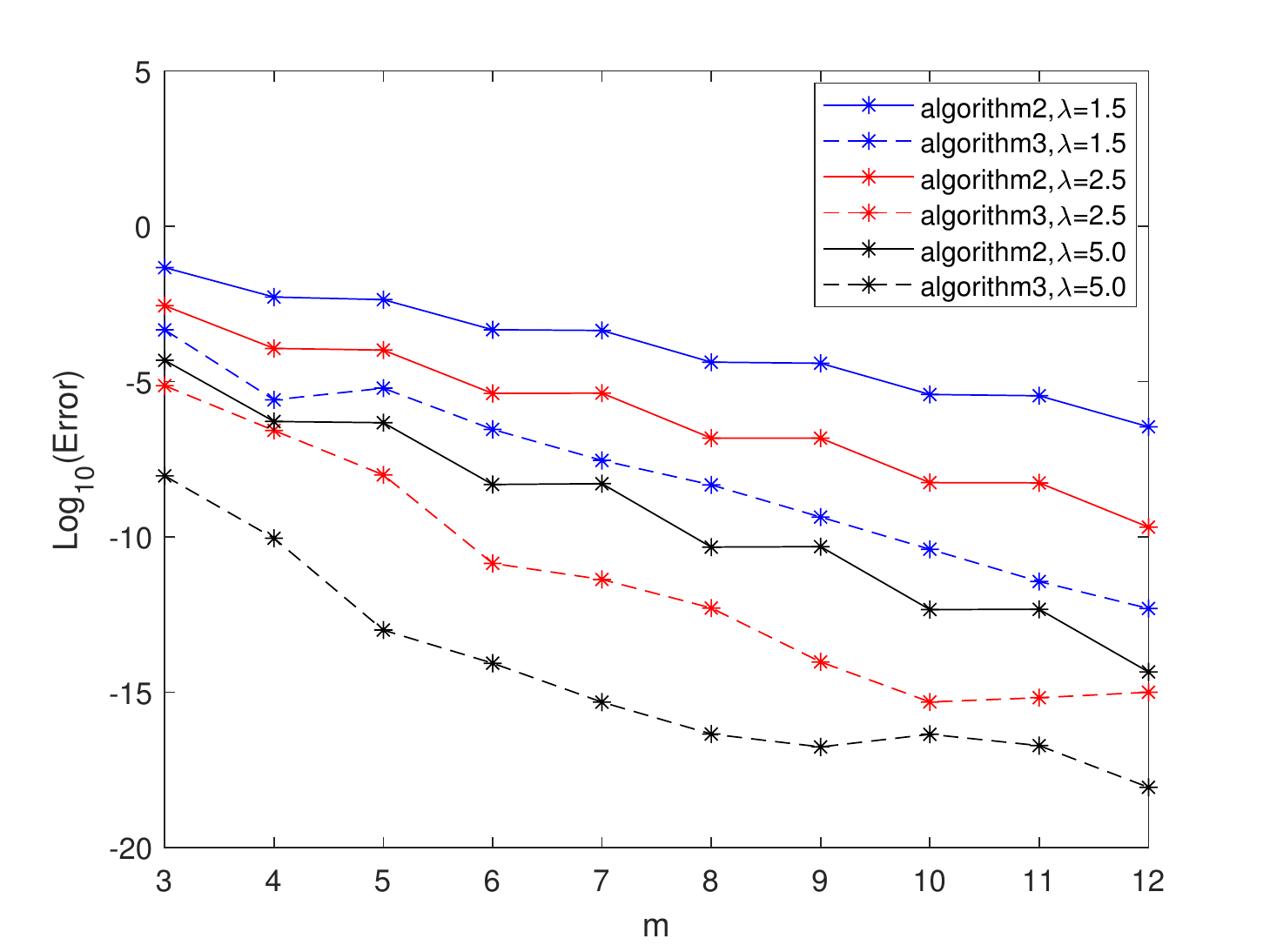}
	\caption{the absolute errors of the approximations to $I_1(p=0,m=7,n=8)$ for different singularities $\xi$} 
	\label{fig:4} 
\end{figure}

We have presented efficient method for the computation of Cauchy principal value integrals and Hadamard finite-part integrals. The error estimation and convergence analysis for the corresponding method is also given. The new scheme based on numerical divided difference is of great accuracy and avoid the cancellation caused by the singular point $\xi$ effectively. Furthermore, the cycle index of the symmetric group allows evaluating the quadrature more efficiently.

\section*{Acknowledgements}

The author would like to thank Prof. Xinghua Wang (Zhejiang University, PRC) for his helpful suggestions
and constructive criticisms which greatly help completing the original manuscript.

\end{document}